\documentclass[12pt]{article}
\usepackage[mathscr]{eucal}
\usepackage{amssymb,latexsym}
\usepackage{verbatim}
\usepackage{amsmath}
\usepackage{amsthm}
\usepackage{pdfsync}

\usepackage{color}

\usepackage[normalem]{ulem}

\usepackage{hyperref}

\usepackage{authblk}

\newtheorem{thm}{Theorem}
\newtheorem{prop}[thm]{Proposition}
\newtheorem{conj}[thm]{Conjecture}

\newcommand\bbR{{\mathbb R}}

\renewcommand\S{\Sigma}
\newcommand\s{\sigma}
\renewcommand\d{\partial}

\newcommand\e{\epsilon}
\renewcommand\b{\beta}
\renewcommand\r{\rho}

\newcommand\ric{{\rm Ric}}

\newcommand\g{\gamma}
\newcommand\8{\infty}
\renewcommand\a{\alpha}

\newcommand\beq{\begin{eqnarray}}
\newcommand\eeq{\end{eqnarray}}
\newcommand\ben{\begin{enumerate}}
\newcommand\een{\end{enumerate}}
\newcommand\bit{\begin{itemize}}
\newcommand\eit{\end{itemize}}

\newcounter{mnotecount}[section]

\author{Gregory J. Galloway\thanks{email: galloway@math.miami.edu}}

\affil{Department of Mathematics
\\University of Miami}

\title{A note on the Lorentzian splitting theorem}

\begin{document}
\date{}
\maketitle  

\begin{abstract}  
We present a version of the Lorentzian splitting theorem under  a weakened Ricci curvature condition.  The proof  makes use of basic properties of achronal limits \cite{horo1, horo2}, together with the geometric maximum principle for $C^0$ spacelike hypersurfaces  in \cite{AGH}.  Our version strengthens a related result in \cite{Yun} in the globally hyperbolic setting by removing a certain boundedness condition on the Ricci curvature.
\end{abstract}

\section{Introduction}

In his well-known problem section from the early 1980's \cite{Yau}, S.-T. Yau posed the problem of establishing  a Lorentzian analogue of the Cheeger-Gromoll splitting theorem.  This problem was completely resolved by the end of the 1980's in a series of papers  \cite{BeemetalLST, Esch, GalLST, Newman}.  The most basic version, due to Eschenburg, \cite{Esch}, is as follows.

\begin{thm} [Lorentzian Splitting Theorem] \label{LST}
Let $M$ be a globally hyperbolic, timelike geodesically complete spacetime, satisfying the strong energy condition, $\mathrm{Ric}(X,X) \ge 0$, for all timelike $X$. If $M$ admits a timelike line, then $M$ splits as a metric product, i.e.
$(M, g)$ is isometric to  $(\bbR \times S, -dt^2 + h)$  where $S$ is a smooth, geodesically complete, spacelike hypersurface with induced metric $h$.
\end{thm}

The Lorentzian splitting problem  posed by Yau was in fact motivated by his idea to establish the {\it rigidity} of the Hawking-Penrose singularity (incompleteness) theorems.  This viewpoint became formalized in a conjecture of R. Bartnik \cite{Bart88}.

 \begin{conj} [Bartnik Splitting Conjecture] \label{Bartnikconj} 
 Suppose that $M$ is a globally hyperbolic spacetime, with compact Cauchy surfaces, which satisfies the strong energy condition. If $M$ is timelike geodesically complete, then 
$M$ splits, i.e. $(M,g)$ is isometric to 
$(\bbR \times S, -dt^2 + h)$, where $S$ is a smooth spacelike Cauchy hypersurface, with induced metric $h$.
\end{conj}  

Thus, according to the conjecture, such spacetimes  can fail to be singular (timelike geodesically incomplete) only under very special circumstances.  To establish the conjecture it is sufficient to establish the existence of a timelike line or a CMC Cauchy surface.  In the former case, one can apply the Lorentzian splitting theorem, whereas in the latter case, one can apply Corollary 1 in \cite{Bart88}.  There are various conditions that ensure either the existence of a timelike line or a CMC Cauchy surface. Nevertheless, the conjecture remains open; for further discussion, see e.g., \cite{GalBart, GalLing24, horo2}. 

Meanwhile, over the years various Hawking-Penrose type singularity results have been proven under weakened (or averaged) energy conditions; see e.g. \cite{Tipler1, Ehrlich, Borde1}, and more recently \cite{Fewster1,Fewster2}, and references therein.  In this note we present a proof of  the Lorentzian splitting theorem under a weakening of the strong energy condition, as discussed in the next section.   
We make use of the properties of {\it achronal limits}~\cite{horo1, horo2}, as they enable certain simplifications and allow much of the proof to remain at a basic causal-theoretic level.
These properties
have been developed for globally hyperbolic spacetimes.  For comments regarding the non-globally hyperbolic case, see the first remark at the end of Section~2.

\section{A Lorentzian splitting theorem}

We briefly recall some basic notions and facts
from the causal theory of spacetime.  See e.g.  \cite{Penrose, HE, ON} 
for further details.  

By a spacetime we mean a smooth time oriented Lorentzian manifold $(M,g)$ of dimension $n+1$, $n \ge 2$. 
For $p \in  M$, $I^+(p)$, the
timelike future of $p$ (resp., $J^+(p)$, the causal future of $p$) is the set of
points $q\in M$ for which there exists a future directed
timelike (resp., future directed causal) curve from $p$ to $q$.
Since small
deformations of timelike curves remain timelike, the sets $I^+(p)$
are  open.
More generally, for $A\subset M$, $I^+(A)$, the timelike  future of $A$,  (resp.,
$J^+(A)$, the causal future of $A$) is the set
of points $q\in M$ for which there exists a future directed
timelike (resp., future directed causal) curve from  a point $p\in A$ to~$q$.
Note that $I^+(A)= \cup_{p\in A}I^+(p)$, and hence the sets $I^+(A)$
are open.  The timelike and causal pasts $I^-(p)$, $J^-(p)$,
$I^-(A)$, $J^-(A)$ are defined in a time dual manner.
An {\it achronal boundary} is a non-empty set of the form $A = \partial I^\pm(S)$, for some subset $S \subset M$. An achronal boundary is in general an edgeless achronal $C^0$ hypersurface of $M$ (\cite{Penrose}, \cite{ON}).

Classically, a spacetime $(M,g)$ is globally hyperbolic provided (i) it is strongly causal (i.e. there are no closed or `almost closed' causal curves) and (ii) all `causal diamonds' $J^+(p) \cap J^-(q)$ are compact. A Cauchy (hyper)surface for a spacetime $(M,g)$ is a set $S$ that is met by every inextendible causal curve exactly once.  A Cauchy surface $S$ is necessarily an edgeless acausal $C^0$ hypersurface. By considering the flow of a timelike vector field, it can be see  that any two Cauchy surfaces are homeomorphic, and that if $S$ is a Cauchy surface for $M$ then $M$ is homeomorphic to $\bbR \times S$.  It is a basic fact that a spacetime $(M,g)$ is globally hyperbolic if and only if it contains a Cauchy surface $S$.   

We will make use of the notion of a future (resp. past) ray. Let $\g:[0,a) \to M$, $a \in (0,\infty]$,  be a future inextendible causal curve.  
$\g$ is a {\it future ray} if each segment of $\g$ is maximal, i.e. each segment realizes the Lorentzian  distance \cite{BEE} between its end points.
Given a subset $S \subset  M$, 
$\g$  is a {\it future $S$-ray} if $\g$ starts on $S$ and each initial segment of $\g$ realizes the Lorentzian distance to $S$. A future ray or $S$-ray is either a timelike or null geodesic (when suitably parameterized).  Past rays and past $S$-rays are defined in an analogous manner.

 \smallskip
 The main aim of this section is to prove the following splitting result.

\begin{thm}\label{split}  Let $(M,g)$ be a globally hyperbolic, timelike geodesically complete spacetime, such that along each (past and future) timelike ray $\a: [0,\infty) \to M$, the Ricci curvature satisfies 
 the following curvature condition:
\beq\label{eq:curv}
\liminf_{t \to \infty} \int_0^{t} \ric(\a'(t),\a'(t))\, dt \ge 0 \, .
\eeq
If $M$ admits a timelike line $\g$, then $M$ splits, i.e.\
$(M, g)$ is isometric to  $(\bbR \times S, -dt^2 + h)$  where $S$ is a smooth, geodesically complete, spacelike hypersurface with induced metric~$h$.
\end{thm}

\noindent
A related result was obtained in \cite{Yun}; see the comments at the end of this section for further details.

As noted in the introduction, our proof of Theorem \ref{split}  makes use of certain results obtained in \cite{horo1, horo2}.  These results are established via basic causal theory, and, in particular, do not involve an analysis of Lorentzian Busemann functions.  

In the discussion that follows it's assumed that spacetime $(M,g)$ is globally hyperbolic.

Let   $\g:[0,\infty) \to M$, be a future complete unit speed timelike ray. 
As defined in~\cite{horo1}, the \emph{past ray horosphere} with respect to 
$\g$, $S^-_\infty$, is given by,
$$
S^-_\infty \,=\, \d \left(\bigcup_{k = 1}^\infty I^-(S^-_k) \right),
$$
where $S^-_k$ is the past Lorentzian sphere,
$$
S^-_k = \{x \mid d\big(x, \g(k) \big) =k\} \,,
$$
where $d$ is the Lorentzian distance function.  $S^-_\infty$ enjoys the following basic properties.  

\begin{prop}[{\cite[Sec. 3]{horo1}}]\label{rays}
$S^-_\infty(\gamma)$ is an edgeless achronal $C^0$ hypersurface which admits a future 
$S^-_\infty(\g)$-ray from each of its points.   In particular,
 $\g$ is itself an $S^-_\infty(\gamma)$-ray. 
\end{prop}

Note that $\bigcup_{k = 0}^\infty I^-(S^-_k)  = I^-\left(\bigcup_{k = 0}^\infty (S^-_k) \right)$.   Hence $S^-_\infty$ is an achronal boundary, and as such  is an edgeless achronal $C^0$ hypersurface.  The future $S$-ray guaranteed at each point can be either a timelike or null geodesic.  In the latter case it is necessarily contained in $S^-_\infty$.   We call $p \in S^-_\infty$ a \emph{null point} if there is a future null $S^-_\infty$-ray from $p$. Otherwise, we call $p$ a \emph{spacelike point} of $S^-_\infty$. 
The set of spacelike points is open and 
acausal  in $S^-_\infty$. 

Future ray horospheres, with respect to past complete timelike rays, are defined in a time dual manner,
$
S^+_\infty \,=\, \d \left(\bigcup_{k = 0}^\infty I^+(S^+_k) \right),
$
and satisfy analogous properties.  Also, the points of a future ray horosphere are classified, as null or spacelike, in a similar  manner.

The next result describes an important structural connection between past and future ray horospheres; see \cite[Lem. 4.10]{horo2}.

\begin{prop} \label{tangent} Suppose that $S^-_\8$ and $S^+_\8$ are past and future ray horospheres, respectively, satisfying $I^+(S^+_\8) \cap I^-(S^-_\8) = \emptyset$. Then for any point $p \in S^-_\8 \cap S^+_\8$, one of the following  holds:

\ben
\item [(1)] $p$ is a spacelike point for both horospheres, and there is a unique future $S^-_\infty$-ray from $p$, and a unique past $S^+_\infty$-ray from $p$, both of which are timelike and join to form a timelike line. 
\item [(2)] $p$ is a null point for both horospheres, and there is a unique future $S^-_\infty$-ray from $p$, and a unique past $S^+_\infty$-ray from $p$, both of which are null and join to form a null line, $\b$, with $\b \subset S^-_\infty \cap S^+_\infty$.
\een
\end{prop}

Propositions \ref{rays} and \ref{tangent} are proved by basic causal theoretic methods. 
The key geometric analytic tool needed for the proof of Theorem \ref{split} involves the notion of mean curvature inequalities {\it in the support sense}; see e.g. \cite{AGH, horo1, horo2}.

Let $S$ be a  locally achronal $C^0$ hypersurface in $M$.
We say that $\Sigma$ is a {\it future support hypersurface} for $S$ at $q \in S$, if $q \in \Sigma$ and $\Sigma$ is a smooth spacelike hypersurface such that, in a neighborhood of $q$, $\Sigma \subset J^+(S)$.  
Further, we say that  $S$ has mean curvature $\le a$ {\it in the support sense} if for all $q \in S$ and $\epsilon > 0$, there is a  future support  hypersurface 
$\S_{q,\epsilon}$ to $S$ at $q$ with mean curvature $H_{q, \epsilon}$ satisfying 
\beq\label{mean}
H_{q, \epsilon}(q) \le a + \epsilon  \,.
\eeq
(By our sign conventions, $H_{q, \epsilon}$ is the divergence along 
$\S_{q,\epsilon}$ of the future directed unit normal field.)
Time-dually, we have the notion of a  past support hypersurface for $S$ at $q\in S$, and can refer to $S$ as having
mean curvature $\ge a$ in the support sense by requiring the reverse inequality, 
$H_{q, \epsilon}(q) \ge a - \epsilon$.

\smallskip
The following theorem  generalizes
aspects of Theorem 4.2  in \cite{horo1}, in particular its extension to the curvature condition  \eqref{eq:curv}.

\begin{thm}\label{convex}
Let $(M,g)$ be a globally hyperbolic, future timelike geodesically complete spacetime, and let  $S$ be a locally achronal $C^0$ hypersurface.  Suppose at each point $p \in S$, there exists a  future complete timelike 
$S$-ray $\a: [0,\infty) \to M$ starting at $\a(0) = p$ such that the curvature condition \eqref{eq:curv} holds. Then $S$ has mean curvature $H \ge 0$ in the support sense.
\end{thm}

\proof Fix $p \in S$ and a future complete unit speed timelike $S$-ray $\a$ with 
$\a(0) = p$.  
Since for any $t$, $\a|_{[0,t]}$ maximizes distance to $S$, one sees that the past distance sphere $S^-_t(\a(t))$  lies locally to the past of $S$ near $q$.
Further, since there are no cut points to $\a(t)$ along $\a|_{[0,t]}$, the past sphere 
$S^-_{t-r}(\a(t))$, $r \in [0,t)$, is smooth near $\a(r)$. It follows, in particular, that
$S^-_t(\a(t))$ is a past support hypersurface for $S$ near $p \in S$.  

\smallskip
The theorem is a consequence of the following claim.

\smallskip
\noindent
{\it Claim:}  For any $\e > 0$ there exists   $T > 0$ such that the distance sphere $S^-_T(\a(T))$ has mean curvature $H \ge -\e$ at $p$.

\proof[Proof of the claim (cf.\ \cite{Fewster1}).] 
Suppose not; i.e.\ suppose there exists $\e > 0$ such that for all arbitrarily large $T$, $H < -\e$ at $p$.  
By the curvature assumption, there exist $r_0 > 0$ such that along $\a$, 
\beq\label{eq:curvineq}
\int_0^r \ric(\a'(\hat{r}),\a'(\hat{r})) d\hat{r} \ge -\e \,, \quad \text{for all } r \ge r_0 \,.
\eeq
We can assume that $T > r_0$.
Let $H(r)$ be the mean curvature of the distance sphere $S^-_{T-r}(\a(T))$ at 
$\a(r)$, $r \in [0,T)$.    By the Raychaudhuri equation (traced Riccati equation \cite{HE}), along $\a|_{[0,T)}$ we have,
\beq\label{eq:raychaud}
\frac{dH}{dr} = - \left(\r + \frac{H^2}{n}\right)   
\eeq 
where $\r = \ric(\a', \a') + \s^2$, and where $\s$, the shear tensor, is the trace free part of the second fundamental of $S^-_{T-r}(\a(T))$  at $\a(r)$.

Since $\s^2 \ge 0$, \eqref{eq:curvineq} implies, 
\beq\label{eq:curvineq2}
\int_0^r \r(\hat{r}) d\hat{r} \ge -\e \,, \quad \text{for all } r \in [r_0,T) \,.
\eeq
Integrating \eqref{eq:raychaud} from $0$ to $r$  and making use of \eqref{eq:curvineq2}, we obtain,
\beq\label{eq:int}
- H(r) =   \int_0^r  \frac{H^2}{n} d\hat{r} +  \int_0^r \r d\hat{r}  - H(0) >  \int_0^r  \frac{H^2}{n} d\hat{r} \, ,
\eeq
for $r \in [r_0,T)$.
Set $Q(r) = \int_0^r \frac{H^2}{n} d\hat{r}$; note that $Q(r) > 0$ on $[r_0, T)$. Then from \eqref{eq:int}
we obtain the differential inequality,
\beq
\frac{dQ}{dr}  = \frac{H^2}{n}  >  \frac{Q^2}{n}  \,\,\, (> 0)
\eeq 
for $r \in [r_0,T)$.   Dividing by $Q^2$ and integrating from $r_0$ to $r$ gives,
\beq
\frac{1}{Q(r_0)}> \frac{1}{Q(r_0)} - \frac{1}{Q(r)} >
\int_{r_0}^r \frac{1}{n} d\hat{r} 
= \frac{r -r_0}{n} \,.
\eeq
By taking $T$ sufficiently large and $r$ sufficiently close  to $T$, the inequality
$$ 
\frac{r -r_0}{n}  <  \frac{1}{Q(r_0)}
$$ will be violated, and so we obtain a  contradiction. This establishes the claim and the proof of the theorem.\qed 

\medskip
We now proceed to the proof of Theorem \ref{split}, which is modeled on the proof of \cite[Theorem 4.1]{horo2}.

\proof[Proof of Theorem \ref{split}]  

Let $\g : (-\8, \8) \to M$ be a complete future-directed unit speed timelike  line. We may decompose $\g$ into a future ray $\g_+:[0,\infty)$, $\g_+(t) = \g(t)$, and a past ray  $\g_-:[0,\infty)$, $\g_-(t) = \g(-t)$.  Associated to these rays are the past ray horosphere $S^-_\8 = S^-_\8(\g_+)$ and the future ray horosphere  $S^+_\8= S^+_\8(\g_-)$, respectively.

 It follows from Proposition \ref{tangent} that $p = \g(0) \in S^-_\8 \cap S^+_\8$ is  a spacelike point for $S^-_\8$ (as well as $S^+_\8$), and hence there exists an acausal neighborhood 
 $U^- \subset S^-_\8$ consisting entirely of spacelike points.  In particular, there exist a future complete timelike $S^-_\8$-ray based at each point of $U^-$.  Theorem \ref{convex} then implies that $U^-$ has mean curvature $H \ge 0$ in the support sense.  In a time dual manner, there exists an acausal neighborhood 
 $U^+ \subset S^+_\8$ having mean curvature $H \le 0$ in the support sense.
 Moreover, as  $I^+(S^+_\8) \cap I^-(S^-_\8) = \emptyset$ 
 (which follows easily from \cite[Proposition 2.5]{horo1}), $U^+$ lies locally to the future of $U^-$ near $p$.  Then by the geometric maximum principle, \cite[Theorem 3.6]{AGH},\footnote{By \cite[Proposition 3.5]{AGH}, a certain technical condition in the statement of  \cite[Theorem 3.6]{AGH} is satisfied in our situation.}, there is a smooth spacelike hypersurface $U \subset U^+ \cap U^-$ having mean curvature $H = 0$. 
 
Consider the normal exponential map $\Phi: \bbR \times U  \to M$,
\beq\label{exp}
\Phi(t, q) = \exp_q tu \,,
\eeq
where $u =$ the future directed unit timelike normal to $U$.  It follows from Proposition~\ref{tangent} that for each $q \in U$, $\a_+(t) = exp_q tu$, $t \in [0,\8)$, is a future timelike $U$-ray, and hence is free of focal points to $U$.  Similarly,   $\a_-(t)  =exp_q(- tu)$, $t \in [0,\8)$, is a past timelike $U$-ray free of focal points. Since, in addition, such $U$-rays can't intersect,  it follows that 
$\Phi: \bbR \times U \to M$ is a diffeomorphism onto its image $W_U \subset M$. 

Now for each fixed $t \in \bbR$, consider the level set 
$U_t = \Phi(\{t\} \times U) =$ 
$\{\exp_q tu : q \in U\}$.  Let $H(t)$ be  the mean curvature of $U_t$. By the Raychaudhuri equation, along each future directed  normal geodesic $\a_+$ starting on $U$, we have,
\beq\label{eq:raychaud2}
\frac{\d H}{\d t} = -\ric(\a_+',\a_+') - \frac{1}{n}H^2 -\s^2 \,. 
\eeq
By an argument similar to the proof of the claim in the proof of 
Theorem~\ref{convex} 
one must have $H(t) \ge 0$ for all $t \in [0,\8)$, as otherwise $H(t)$ could not be defined for all 
$t \in [0,\8)$.   

Suppose $H(t_0) > 0$ for some $t_0 \ge 0$.  Let $\tilde H(t) = - H(t)$.  Since $\tilde H(t_0) < 0$, essentially the same argument shows that $\tilde H(t)$, and hence $H(t)$, can't be defined for all $t \in (-\8, t_0)$.  It follows that $H(t) = 0$ for all $t \in [0,\infty)$.  Then \eqref{eq:raychaud2} gives,
$$
\ric(\a_+',\a_+') = -\s^2 \le  0 \,, \quad \text{for all } t \in [0,\8) \,.  
$$
Taken together with the curvature condition  \eqref{eq:curv}, we conclude that $\ric(\a_+',\a_+') = -\s^2 = 0$ for all  $t \in [0,\8)$.   This implies that $U_t$ is  totally geodesic for  each $t \in [0,\8)$.  Time dual arguments also show that $U_t$ is totally geodesic for $t \in (-\infty,0]$.   Standard arguments now imply that 
$(W_U, g)$ is isometric to the product $(\bbR \times U, -dt^2 +h)$, where $h$ is the induced metric on $U$.  

We now want to extend this to a global splitting.  Let $V \supset U$ be the connected component of the set of common spacelike points in  
$S^-_\8 \cap S^+_\8$ containing $p$; $V$ is acausal and open in  $S^-_\8 \cap S^+_\8$.   From Proposition \ref{tangent}, based at each point $q \in V$  is a unique future $S^-_\infty$-ray $\b_+$ and a unique past $S^+_\infty$-ray $\b_-$, which join together to form a timelike line $\b$. Then, by the same arguments as above, we obtain a local spacetime splitting around the timelike line $\b$ similar to that above. 
It follows that $V$ is a smooth totally geodesic spacelike hypersurface and that 
$(W_V = \Phi(\bbR \times V) , g)$ is isometric, via the normal exponential map $\Phi$ extended to $V$, to the product $(\bbR \times V, -dt^2 +h)$, where now $h$ is the induced metric on $V$.
From the fact that $V$ is totally geodesic and that $M$ is timelike geodesically complete, it is shown in the proof of 
 \cite[Theorem~4.11]{horo2}  that $V$ is geodesically complete. We include this argument in the appendix.
 It then follows that  $(\bbR \times V, -dt^2 +h)$ is geodesically complete and that $\{0\} \times V$ is a Cauchy surface in $(\bbR \times V, -dt^2 +h)$
 (see  \cite[Theorems 3.67 and 3.69]{BEE}).  Hence,
$(W_V,g)$ is geodesically complete and, moreover,  $W_V =  D(V)$, where $D(V)$ is the domain of dependence of $V$ within $M$.  One then easily shows  that $V$ is a Cauchy surface in $M$ (e.g., by showing that there is no Cauchy horizon, $H^+(V) = H^-(V) = \emptyset$), and hence $W_V = M$.  Thus, 
$(M ,g)$ is isometric to  $(\bbR \times V, -dt^2 +h)$.\qed
 
\smallskip

We mention one consequence of Theorem \ref{split}.  Let $(M,g)$ be a future timelike geodesically complete spacetime with compact Cauchy surfaces.  By the time dual of \cite[Prop. 5.9]{horo2}, if the future
causal boundary of $M$ is spacelike\footnote{For a discussion of the causal boundary of spacetime see e.g., \cite[Sect. 6.8]{HE}.  The future causal boundary is spacelike if no TIP is properly contained in another TIP.} then $M$ admits a timelike line.  Hence, under this additional assumption,  Conjecture \ref{Bartnikconj} holds, as one can apply the standard Lorentzian splitting theorem.  In a similar manner, one has the following generalization. 

 \begin{thm}\label{bartsplit} 
 Suppose that $(M,g)$ is a timelike geodesically complete spacetime with compact Cauchy surfaces, which satisfies the energy condition \eqref{eq:curv}.  If the future causal boundary is spacelike then $M$ splits, i.e. $(M,g)$ is isometric to 
 $(\bbR \times S, -dt^2 + h)$, where $S$ is a smooth spacelike Cauchy hypersurface, with induced metric $h$. 
\end{thm}

 Now it is Theorem \ref{split} that insures that $(M,g)$ splits.  For further relevance of the future causal boundary being spacelike see \cite[Sec. 5.1]{GalLing24}.

\medskip
\noindent
{\it Final remarks.}

\smallskip
\noindent
1.  The causal theoretic techniques employed in the proof of Theorem \ref{split} circumvent a number of the technical issues involving the regularity of Lorentzian Busemann functions. Working instead with `ray horospheres'  leads to a number of simplifications.  However, using the  regularity theory of Lorentzian Busemann functions for timelike geodesically complete spacetimes, as developed in \cite{GH}, one can obtain a proof of Theorem~\ref{split} without the global hyperbolicity assumption.  Such an approach was taken in \cite{Yun}, but the arguments there required an additional condition on the Ricci curvature. Roughly speaking, it is required that the integral along timelike rays of the negative part of the Ricci curvature be bounded above.   By combining results from \cite{GH}, together with arguments in this paper, this extra condition can be avoided.  We make a brief comment about this.

From results in \cite{GH} there exists a neighborhood $U$ of $\g(0)$ ($\g$ being the given timelike line) such that the local Busemann level sets 
$S_{\pm} = \{b^{\pm} = 0 \} \cap U$ are acausal $C^0$ spacelike hypersurfaces passing through 
$\g(0)$, with $S_-$ lying to the causal future  of $S_+$.  Using the local regularity of certain Busemann support functions as discussed in~\cite{GH} (see also \cite{BEE}), together with arguments similar to those in the proof of Theorem~\ref{convex}, one can show that $S_-$ has mean curvature $\le 0$ in the support sense
and $S_+$ has mean curvature $\ge 0$ in the support sense.  Then, applying the geometric maximum principle in \cite{AGH}, as in the proof of Theorem \ref{split}, one has that $S_-$ and $S_+$ agree along a smooth maximal spacelike  hypersurface $S$ in $U$ 
passing through $\g(0)$.  
Now, making use of the existence of past and future $S$-rays at each point of $S$ (see \cite{GH, BEE}), one may proceed as in the proof of Theorem \ref{split} to obtain a local splitting.  This local splitting can be extended to a global one as in \cite{Esch} or \cite{BEE}.

\smallskip
\noindent
2. In \cite{GalCG}, a Riemannian analogue of Theorem \ref{split} was obtained by modifying the proof of Eschenburg and Heintze \cite{EH} of the Cheeger-Gromoll splitting theorem.   This suggests the possibility that the recent elliptic methods \cite{mcc} based on the $p$-Laplacian in proving the Lorentzian splitting theorem may be relevant to the situation considered here. 

\smallskip
\noindent
3.  A quite different direction for generalizing the Lorentzian splitting theorem lies in the current, very active area of research on low-regularity Lorentzian geometry. This includes cases where the Lorentzian metric has lower than $C^2$ regularity, as in low regularity versions of the Hawking-Penrose singularity theorems (see \cite{Graf} and references therein), or where synthetic approaches might be used \cite{Kunz}. See, for example, \cite{Solis}, which  extends the Lorentzian splitting theorem in the sectional curvature case to a synthetic setting, and \cite{mcc} which describes how the novel methods developed there will lead to low-regularity versions of the Lorentzian splitting theorem under Ricci curvature bounds.  In the present context, one might  wonder whether the causal-theoretic techniques developed in \cite{horo1,horo2}, together with metric `regularization', e.g.,  as in \cite{KunzStein}, could  possibly play a role in obtaining a low-regularity version of the Lorentzian splitting theorem.

 \section*{Appendix}
 
 Let the setting and notation be as in  the statement and proof of Theorem \ref{split}, 
and let $V \subset S^-_\8 \cap S^+_\8$ be as defined in the proof.  
We want to show that $V$ is geodesicallly complete.  Fix $p \in V$, and choose any $X \in T_pV$ with $|X| = 1$.  Let $\s: [0,a) \to M$, $a \in (0,\8]$, be the unique unit speed geodesic with  $\s(0) = p$ and $\s'(0) = X$, which is maximally extended in the direction $X$ in $M$.  We first observe that $\s$ does not leave $V$. 
 Because $V$ is totally geodesic, $\s$ is initially contained in $V$. Fix any $s_0  \in (0,a)$ with 
 $\s([0, s_0)) \subset V$. Again, since $V$ is totally geodesic, the future unit normal field $u$ of $V$ is parallel along $\s|_{[0,s_0)}$. 
 By Proposition \ref{tangent}, there is a unique future timelike $S^-_\infty$-ray 
 $\g_x$ from each $x \in V$, with $\g_x'(0) = u_x$.  
 Whether or not $q = \s(s_0)$ lies in $V$, there is a well-defined limit vector $u_q = \lim_{s \to s_0}u_{\s(s)}$, obtained by parallel transporting $u$ on all of $\s_{[0,s_0]}$, with $u_q$ necessarily future timelike. Let $\g_q$ be the future-directed unit speed timelike geodesic with 
 $\g_q'(0) = u_q$, which is necessarily complete.    Since $q \in \overline{V} \subset S^-_\8 \cap S^+_\8$, $\g_q$ is a curve from $S^-_\8$. Moreover, since  $\g_{\s(s)}$ is an $S^-_\8$-ray for all $s \in [0, s_0)$, it follows that $\g_q|_{[0, \8)}$ is an $S^-_\8$-ray, as well. Since $q \in S^-_\8 \cap S^+_\8$, Proposition \ref{tangent} implies that $q$ is a spacelike point, and hence $q \in V$.  This shows that $\s$ can never leave $V$, i.e., we have $\s : [0, a) \to V$. 
 
 Now we show  that $\s$ is complete in the direction $X$, i.e.\ that $a = \infty$.
 Suppose to the contrary that $a < \8$. Then the curve $c(s) = \Phi(-2s, \s(s))$, 
$c : [0, b) \to W_V \subset M$, 
 is a past-directed timelike geodesic in $M$, and 
 $\s(s) = \mathrm{exp}_{c(s)}(2s \Phi_*(\d_t))$. By timelike completeness, $c$ extends to $[0,a]$. Furthermore, the vector field $\Phi_*(\d_t)$ is parallel in $W_V$, and hence, by parallel translating along $c$, has a limit at $c(a)$. Hence, 
$\s(s) = \mathrm{exp}_{c(s)}(2s \Phi_*(\d_t))$
 has a limit as $s \to a$, i.e., $\s$ extends continuously, and hence as a geodesic to $[0, a]$, contrary to the definition of $a$.
 Hence,  $a = \8$.  Since $X$ was an arbitrary unit vector at $p$ in $V$, we have, in fact, shown that 
 $\exp_p: T_pV \to V$ is defined on all of $T_pV$.   Thus by Hopf-Rinow, $V$ is geodescially complete.
 
\bigskip
 
\noindent
\textsc{Acknowledgements.}  We are grateful to Carlos Vega for many helpful comments on an earlier version of this paper.

\providecommand{\bysame}{\leavevmode\hbox to3em{\hrulefill}\thinspace}
\providecommand{\MR}{\relax\ifhmode\unskip\space\fi MR }
\providecommand{\MRhref}[2]{%
  \href{http://www.ams.org/mathscinet-getitem?mr=#1}{#2}
}
\providecommand{\href}[2]{#2}

\bibliographystyle{amsplain}
\bibliography{split}

\end{document}